\documentclass[11pt,a4paper]{article}

\usepackage{graphics,epsfig,theorem,latexsym,amssymb,amsmath, amsfonts}
\usepackage{times,epic,eepic}

\usepackage{pifont}

\usepackage{cite}
\usepackage{url}

\usepackage{color}

\usepackage{newlfont}

 \newenvironment{keywords}{ \noindent {\small\bf Key Words}:}{ }

\def\bd{\begin{description}}
\def\ed{\end{description}}

\def\beq{\begin{equation}}
\def\eeq{\end{equation}}

 \newtheorem{theorem}{Theorem}

 \newtheorem{example}{Example}

\begin{document}
%\begin{article}

\title{Higher order numerical differentiation\\  on the Infinity Computer}

%% Author name
%\newcommand{\nms}{\normalsize}
\author{Yaroslav D. Sergeyev\footnote{Yaroslav D.
Sergeyev, Ph.D., D.Sc., is Distinguished Professor at the University
of Calabria, Rende, Italy.
 He is also Full Professor (a part-time contract) at the N.I.~Lobatchevsky State University,
  Nizhni Novgorod, Russia and Affiliated Researcher at the Institute of High Performance
  Computing and Networking of the National Research Council of Italy.   }}

\date{}

\maketitle

 \begin{abstract}
There exist many   applications where  it is necessary to
approximate numerically derivatives of a function   which is given
by a computer procedure. In particular, all the fields of
optimization have a special interest in such a kind of
information. In this paper, a new way to do this is presented for
a new kind of a computer -- the Infinity Computer -- able to work
numerically with finite, infinite, and infinitesimal numbers. It
is proved that the Infinity Computer is able to calculate values
of derivatives of a higher order for a wide class of functions
represented by computer procedures. It is shown that the ability
to compute derivatives of arbitrary order automatically and
accurate to working precision is an intrinsic property  of the
Infinity Computer related to  its way of functioning. Numerical
examples illustrating the new concepts and numerical tools are
given.
 \end{abstract}


%\subclass{65K05, 90C26, 90C56}

\section{Introduction}
\label{s_m1}

In many practical applications related to the scientific computing
(e.g., in global and local optimization, numerical simulation,
approximation,  etc.)  it is necessary to calculate
  derivatives of a function $g(x)$ which is given by a computer
procedure calculating its   approximation $f(x)$. Very often a
user working with the computing code $f(x)$ is not the person who
has written this code. As a result, for the user  the program
calculating $y=f(x)$ is just a black box, i.e., if it has as the
input a value $x$ then the program returns the corresponding value
$y$ and the user does not know the internal structure of the
program. As a result, when for solving an applied problem the
usage of derivatives is required and a procedure for evaluating
the exact value of $f'(x)$ is not available, we face the necessity
to approximate $f'(x)$ in a way.

In particular, this situation happens very often in the black box
global and local optimization (see
\cite{Corliss_et_al,Pardalos&Resende(2002),fizmatlit,Strongin_Sergeyev})
and related application areas. Let us give a simple but important
example (see, e.g., \cite{fizmatlit,Strongin_Sergeyev,Wolfe})
related to   the problem of finding  the minimal root of an
equation $f(x)=0$ where $x \in [a,b]$ and $f(x)$ is multiextremal
(as a result, there can be several different roots over $[a,b]$),
given by a computer program and such that $f(a)>0$.
  This problem arises in many applications,
such as   time domain analysis (see \cite{kn:circuits}), filter
theory (see \cite{kn:filters}), and wavelet theory (see
\cite{Walnut}) and can be interpreted, for instance, as follows.

It is necessary to know the behavior of a device over a time
interval $[a,b]$. The device starts to work at the time $x=a$ and
it functions correctly while for $x \ge a$ the computer procedure
calculating $f(x)$ returns values $f(x)>0$. Of course, at the
initial moment, $x=a$, the device works correctly and $f(a)>0$. It
is necessary either to find an interval $[a,x^*)$ such that
\begin{equation}
f(x^*)=0,\hspace{.1in} f(x)>0, \hspace{.1in}x \in [a,
x^*),\hspace{.1in} x^* \in (a,b], \label{such that}
\end{equation}
or to prove that $x^*$ satisfying (\ref{such that}) does not exist
in $[a,b]$.  Efficient methods proposed recently for solving this
problem    (see \cite{kn:zero-crossing,kn:SIAM,fizmatlit})
strongly use ideas developed in the field of global optimization.
They require calculating the first derivative $f'(x)$ of $f(x)$
and since a program calculating $f'(x)$ is usually not available,
the problem of finding an approximation of $f'(x)$ arises.

There exist several approaches  to tackle this problem. First,
numerical approximations are used for this purpose (see e.g.,
\cite{Moin} for a detailed discussion). In applications, the
following three simple formulae (more complex and numerically more
expensive approximations can be found in \cite{Moin}) are often
used
 \beq
f'(x)\approx \frac{f(x+h)-f(x)}{h}, \hspace{5mm} f'(x)\approx
\frac{f(x)-f(x-h)}{h},
 \label{m0}
       \eeq
 \beq
f'(x)\approx \frac{f(x+h)-f(x-h)}{2h} \label{m0.1}
       \eeq
by practitioners. However, these procedures are fraught with
danger (see \cite{Moin}) since eventually round-off errors will
dominate calculation. As $h$ tends to zero, both $f(x+h)$ and $
f(x-h)$ tend to $f(x)$, so that their difference tends to the
difference of two almost equal quantities and thus contains fewer
and fewer significant digits. Thus, it is meaningless to carry out
these computations beyond a certain threshold value of $h$.
Calculations of higher derivatives suffer from the same problems.

The complex step method (see \cite{Lyness_Moler}) allows one to
improve  approximations of $f'(x)$ avoiding subtractive
cancellation errors present in (\ref{m0}), (\ref{m0.1}) by using
the following formula to approximate $f'(x)$
 \beq
f'(x) \approx \frac{Im[f(x+ih)]}{h},
 \label{m0.2}
       \eeq
 where $Im(u)$ is the imaginary part of $u$. Though this estimate does not
involve the dangerous difference operation,   it is still an
approximation of $f'(x)$ because it depends on the choice of the
step $h$.

Another approach consists of the usage of symbolic (algebraic)
computations (see, e.g., \cite{Cohen_symbolic}) where $f(x)$ is
differentiated as an expression in symbolic form in contrast to
manipulating   of   numerical quantities represented by the
symbols used to express  $f(x)$. Unfortunately this approach can
be too slow when it is applied to long   codes coming from real
world applications.

There exist an extensive literature (see, e.g.,
\cite{Berz_dual,Bischof,Corliss_et_al} and references given
therein) dedicated to automatic (algorithmic) differentiation (AD)
that is a set of techniques based on the mechanical application of
the chain rule to obtain derivatives of a function given as a
computer program.   By applying the chain rule of derivation to
elementary operations  this approach allows one to compute
derivatives of arbitrary order automatically with the precision of
the code representing $f(x)$.

Implementations of AD can be broadly classified into two
categories that have their advantages and disadvantages (see
\cite{Bischof,Corliss_et_al} for a detailed discussion): (i) AD
tools based on source-to-source transformation changing the
semantics by explicitly rewriting the code; (ii) AD tools based on
operator overloading using the fact that modern programming
languages offer the possibility to redefine the semantics of
elementary operators. In particular, the dual numbers extending
the real numbers by adjoining one new element $d$ with the
property $d^2 = 0$ (i.e., $d$ is nilpotent) can be used for this
purpose (see, e.g., \cite{Berz_dual}). Every dual number has the
form $v = a + db$, where $a$ and $b$ are real numbers and $v$ can
be represented as the ordered pair $(a,b)$. On the one hand, dual
numbers have a clear similarity with complex numbers $z = a + ib$
where $i^2 = -1$. On the other hand, speaking informally it can be
said that the imaginary unit $d$ of dual numbers is a close
relative to infinitesimals (we mean here a general non formalized
idea  about infinitesimals) since the square (or any higher power)
of $d$ is exactly zero and the square of an infinitesimal is
`almost zero'.

All the methods described above use traditional computers as
computational devices and propose a number of techniques to
calculate derivatives on them. In this paper, a new way to
calculate derivatives numerically is proposed. It is made by using
a new kind of a computer -- the Infinity Computer -- introduced in
\cite{Sergeyev_patent,www,informatica} and able to work
\textit{numerically} with finite, infinite, and infinitesimal
quantities. This computer is based on   a new applied point of
view on infinite and infinitesimal numbers (that is not related to
non-standard analysis)   introduced in
\cite{Sergeyev,informatica}. The new approach does not use
Cantor's ideas and works with infinite and infinitesimal numbers
being in accordance with Aristotle's principle `The part is less
than the whole'.

We conclude this introduction by emphasizing that traditional
approaches for differentiation considered above have been
developed \textit{ad hoc} for solving this problem as additional
tools   that should be used together with the traditional
computers. Without these additional tools the traditional
computers are not able to calculate derivatives of functions
defined by computer procedures. In this paper, it is shown that
the ability to compute derivatives of arbitrary order
automatically and accurate to working precision is an intrinsic
property  of the Infinity Computer related to \textit{its way of
functioning}. This is just one of the particular features offered
to the user by the Infinity Computer. Naturally, this is a direct
consequence of the fact that it can execute numerical computations
with infinite and infinitesimal quantities explicitly.

\section{Representation of numbers at the Infinity Computer}
\label{s_m2}

In \cite{Sergeyev,informatica,Korea,Dif_Calculus}, a new powerful
numeral system has been developed to express finite, infinite, and
infinitesimal numbers in a unique framework. The main idea
consists of measuring infinite and infinitesimal quantities   by
different (infinite, finite, and infinitesimal) units of measure.
In this section we give just  a brief introduction to the new
methodology    that can be found in a rather comprehensive form in
the survey \cite{informatica} or in the monograph \cite{Sergeyev}
written in a popular manner.

A new infinite unit of measure   has been introduced    as the
number of elements of the set $\mathbb{N}$ of natural numbers. It
is expressed by a new numeral \ding{172} called \textit{grossone}.
It is necessary to emphasize immediately that the infinite number
\ding{172} is not either Cantor's $\aleph_0$ or $\omega$ and the
new approach is not related to the non-standard analysis. For
instance, one of the important differences consists of the fact
that infinite integer positive numbers that can be viewed by using
numerals including grossone can be interpreted in the terms of the
number of elements of certain infinite sets. Another difference
consists of the fact that \ding{172} has both cardinal and ordinal
properties as usual finite natural numbers.

Formally, grossone is introduced as a new number by describing its
properties postulated by the \textit{Infinite Unit Axiom} (IUA)
(see \cite{Sergeyev,informatica}). This axiom is added to axioms
for real numbers similarly to addition of the axiom determining
zero to axioms of natural numbers when integer numbers are
introduced. Inasmuch as it has been postulated that grossone is a
number,  all other axioms for numbers hold for it, too.
Particularly, associative and commutative properties of
multiplication and addition, distributive property of
multiplication over addition, existence of   inverse  elements
with respect to addition and multiplication hold for grossone as
for finite numbers. This means that  the following relations hold
for grossone, as for any other number
 \beq
 0 \cdot \mbox{\ding{172}} =
\mbox{\ding{172}} \cdot 0 = 0, \hspace{3mm}
\mbox{\ding{172}}-\mbox{\ding{172}}= 0,\hspace{3mm}
\frac{\mbox{\ding{172}}}{\mbox{\ding{172}}}=1, \hspace{3mm}
\mbox{\ding{172}}^0=1, \hspace{3mm}
1^{\mbox{\tiny{\ding{172}}}}=1, \hspace{3mm}
0^{\mbox{\tiny{\ding{172}}}}=0.
 \label{3.2.1}
       \eeq

To express infinite and infinitesimal numbers  at the Infinity
Computer, records similar to traditional positional numeral
systems can be used (see \cite{Sergeyev,www,informatica}). Numbers
expressed in this new positional systems with the radix \ding{172}
are called hereinafter \textit{grossnumbers}. In order to
construct a number $C$ in this system, we subdivide $C$ into
groups corresponding to powers of grossone:
 \beq
  C = c_{p_{m}}
\mbox{\ding{172}}^{p_{m}} +  \ldots + c_{p_{1}}
\mbox{\ding{172}}^{p_{1}} +c_{p_{0}} \mbox{\ding{172}}^{p_{0}} +
c_{p_{-1}} \mbox{\ding{172}}^{p_{-1}}   + \ldots   + c_{p_{-k}}
 \mbox{\ding{172}}^{p_{-k}}.
\label{3.12}
       \eeq
 Then, the record
 \beq
  C = c_{p_{m}}
\mbox{\ding{172}}^{p_{m}}    \ldots   c_{p_{1}}
\mbox{\ding{172}}^{p_{1}} c_{p_{0}} \mbox{\ding{172}}^{p_{0}}
c_{p_{-1}} \mbox{\ding{172}}^{p_{-1}}     \ldots c_{p_{-k}}
 \mbox{\ding{172}}^{p_{-k}}
 \label{3.13}
       \eeq
represents  the number $C$, where finite numbers $c_i\neq0$ called
\textit{grossdigits}   can be   positive or negative. They show
how many corresponding units should be added or subtracted in
order to form the number $C$. Grossdigits can be expressed by
several symbols.

Numbers $p_i$ in (\ref{3.13}) called \textit{grosspowers}  can be
finite, infinite, and infinitesimal, they   are sorted in the
decreasing order with $ p_0=0$
\[
p_{m} >  p_{m-1}  > \ldots    > p_{1} > p_0 > p_{-1}  > \ldots
p_{-(k-1)}  >   p_{-k}.
 \]
 In the record (\ref{3.13}), we write
$\mbox{\ding{172}}^{p_{i}}$ explicitly because in the new numeral
positional system  the number   $i$ in general is not equal to the
grosspower $p_{i}$ (see \cite{informatica} for a detailed
discussion).

 Finite numbers  in this new numeral system are represented
by numerals having only one grosspower $ p_0=0$. In fact, if we
have a number $C$ such that $m=k=$~0 in representation
(\ref{3.13}), then due to (\ref{3.2.1}),   we have $C=c_0
\mbox{\ding{172}}^0=c_0$. Thus, the number $C$ in this case does
not contain grossone and is equal to the grossdigit $c_0$ being a
conventional finite number     expressed in a traditional finite
numeral system.

 The simplest infinitesimal numbers    are represented by numerals $C$
having only finite or infinite negative grosspowers, e.g.,
$6.73\mbox{\ding{172}}^{-4.7}56.7\mbox{\ding{172}}^{-150}$. The
simplest infinitesimal number is
$\frac{1}{\mbox{\ding{172}}}=\mbox{\ding{172}}^{-1}$ being the
inverse element with respect to multiplication for \ding{172}:
 \beq
\mbox{\ding{172}}^{-1}\cdot\mbox{\ding{172}}=\mbox{\ding{172}}\cdot\mbox{\ding{172}}^{-1}=1.
 \label{3.15.1}
       \eeq
Note that all infinitesimals are not equal to zero. Particularly,
$\frac{1}{\mbox{\ding{172}}}>0$ because it is a result of division
of two positive numbers.

In the context of the numerical differentiation discussed in this
paper, it is  worth mentioning that (without going in a detailed
and rather philosophical  discussion on the topic `Can or cannot
dual numbers be viewed as a kind of infinitesimals?') there exist
two formal differences between infinitesimals $C$ from
(\ref{3.13}) and dual numbers (see, e.g., \cite{Berz_dual}).
First, for any infinitesimal $C$ it follows $C^2>0$ (for instance,
$(\mbox{\ding{172}}^{-1})^2>0$) whereas for dual numbers we have
$d^2=0$. Second,   in the context of \cite{Berz_dual} the element
$d$ represented as $(0,1)$ has not its inverse and infinitesimals
$C$ have their inverse.

The simplest infinite numbers     are expressed by numerals having
  positive finite or infinite gross\-powers. They have
infinite parts and can also have  a finite part and infinitesimal
ones.  For instance, the number
\[
 1.5\mbox{\ding{172}}^{14.2}(-10.645)\mbox{\ding{172}}^{5}7.89\mbox{\ding{172}}^{0}81\mbox{\ding{172}}^{-4.2}72.8\mbox{\ding{172}}^{-60}
\]
has two infinite parts $1.5\mbox{\ding{172}}^{14.2}$ and
$-10.645\mbox{\ding{172}}^{5}$ one finite part
$7.89\mbox{\ding{172}}^{0}$ and two infinitesimal parts
$81\mbox{\ding{172}}^{-4.2}$ and $72.8\mbox{\ding{172}}^{-60}$.
All of the numbers introduce above can be grosspowers, as well,
giving so a possibility to have various combinations of quantities
and to construct  terms having a more complex structure.

A working software simulator of the Infinity Computer has been
implemented and the first application -- the Infinity Calculator
-- has been realized. We conclude this section by emphasizing the
following important issue: the Infinity Computer works with
infinite, finite, and infinitesimal numbers \textit{numerically},
not symbolically (see \cite{Sergeyev_patent}).

\section{Numerical differentiation}
\label{s_m3}

Let us return   to the problem of numerical differentiation of a
function $g(x)$. We suppose that   a set of elementary functions
($\sin(x), \cos(x), a^{x}$ etc.) is represented at the Infinity
Computer  by one of the usual ways used in traditional computers
(see, e.g. \cite{Muller})   involving   the argument $x$, finite
constants, and four arithmetical operations. A programmer writes a
program   $P$ that should calculate  $g(x)$ using the said
implementations of elementary functions, the argument $x$, and
finite constants connected by four arithmetical operations.
Obviously, $P$   calculates a numerical approximation $f(x)$  of
the function $g(x)$. As a rule, the programmer does not use
analytical formulae of  $f'(x), f''(x), \ldots f^{(k)}(x)$ to
write the program calculating $f(x)$. We suppose that $f(x)$
approximates $g(x)$ sufficiently well with respect to some
criteria  and we shall not discuss the goodness of this
approximation in this paper.

Then, as often happens in the scientific computing, a user takes
the program $P$ calculating $f(x)$  and is interested to calculate
$f'(x)$ and higher derivatives numerically by using this program.
Computer programs for calculating $f'(x),$ $ f''(x), \ldots $
$f^{(k)}(x)$ and their analytical formulae  are unavailable and
the internal structure of the program calculating $f(x)$ is
unknown to the user.

In this situation, our attention will be attracted to the problem
of a numerical calculation  of the derivatives $f'(x), f''(x),
\ldots f^{(k)}(x)$ and to the information that can be obtained
from the computer procedure $P$ calculating $f(x)$ for this
purpose when it is executed at the Infinity Computer. The
following theorem holds.

\begin{theorem}
\label{t_m1} Suppose that: (i) for a function $f(x)$   calculated
by a procedure implemented at the Infinity Computer  there exists
an  unknown  Taylor expansion in a finite neighborhood $\delta(y)$
of a finite point $y$; (ii) $f(x),$ $f'(x), f''(x), \ldots
f^{(k)}(x)$ assume finite values or are equal to zero for $x \in
\delta(y)$; (iii) $f(x)$ has been evaluated at a point
$y+\mbox{\ding{172}}^{-1} \in \delta(y)$. Then the Infinity
Computer returns the result of this evaluation in the positional
numeral system with the infinite radix \ding{172} in the following
form
 \beq
  f(y+\mbox{\ding{172}}^{-1}) = c_{0}
\mbox{\ding{172}}^{0}   c_{-1} \mbox{\ding{172}}^{-1}  c_{-2}
\mbox{\ding{172}}^{-2}   \ldots c_{-(k-1)}
 \mbox{\ding{172}}^{-(k-1)} c_{-k}
 \mbox{\ding{172}}^{-k},
\label{m1}
       \eeq
      where
 \beq
%\begin{array}{l}
%  f(y) = c_{0},\\
%f'(y) =    c_{-1}, \\
%f''(y)= 2!   c_{-2},\\
% \ldots \\
% f^{(k)}(y)=k!  c_{-k}.
%   \\
%    \end{array}
   f(y) = c_{0}, \,\,
f'(y) =    c_{-1}, \,\, f''(y)= 2!   c_{-2}, \,\,
 \ldots \,\,
 f^{(k)}(y)=k!  c_{-k}.
       \label{m2}
       \eeq
\end{theorem}

\textbf{Proof.} Due to its rules of   operation (see (\ref{3.12}),
(\ref{3.13})), the Infinity Computer   collects  different exponents
of \ding{172} in independent groups $c_{p_{-i}}
\mbox{\ding{172}}^{p_{-i}}$ with finite grossdigits  $c_{p_{-i}}$
when it calculates $f(y+\mbox{\ding{172}}^{-1})$. Since functions
$f(x),$ $f'(x), f''(x), \ldots$ $ f^{(k)}(x)$ assume finite values
or are equal to zero in $\delta(y)$ which is also finite, the
highest grosspower in the number (\ref{m1}) is necessary less or
equal to zero. Thus, the number that the Infinity Computer returns
can have only a finite and infinitesimal parts.

  Four arithmetical operations (see
\cite{Sergeyev_patent,informatica}) executed by the Infinity
Computer with the operands having finite integer grosspowers in
the form (\ref{3.13})  produce only results with finite integer
grosspowers. This fact  ensures that the result
$f(y+\mbox{\ding{172}}^{-1})$ can have only integer non-positive
grosspowers in (\ref{m1}). Due to the rules of the positional
system (see (\ref{3.12}), (\ref{3.13})), the number
$f(y+\mbox{\ding{172}}^{-1})$ from (\ref{m1}) can be written as
follows
\[
  f(y+\mbox{\ding{172}}^{-1}) = c_{0}
\mbox{\ding{172}}^{0}   c_{-1} \mbox{\ding{172}}^{-1}  c_{-2}
\mbox{\ding{172}}^{-2}   \ldots c_{-(k-1)}
 \mbox{\ding{172}}^{-(k-1)} c_{-k}
 \mbox{\ding{172}}^{-k}=
\]
 \beq
  c_{0}
\mbox{\ding{172}}^{0} +  c_{-1} \mbox{\ding{172}}^{-1} + c_{-2}
\mbox{\ding{172}}^{-2} +  \ldots + c_{-(k-1)}
 \mbox{\ding{172}}^{-(k-1)} + c_{-k}
 \mbox{\ding{172}}^{-k}.
\label{m1.1}
       \eeq
The Infinity Computer while calculates the value
$f(y+\mbox{\ding{172}}^{-1})$ does not use the Taylor expansion
for $f(x)$, it just executes commands of the program. However,
this unknown Taylor expansion for $f(x)$ (we emphasize that it is
unknown for: the Infinity Computer itself, for the programmer, and
for the user)  exists in the neighborhood $\delta(y)$ of the point
$y$, for a point $x=y+h \in \delta(y),$ $ h>0$. Thus, it should be
true
 \beq
f(y+h)=f(y)+f'(y)h+ f''(y)\frac{h^2}{2}+ \ldots +
f^{(k)}(y)\frac{h^k}{k!} + \ldots \label{m1.11}
       \eeq
By assuming   $h=\mbox{\ding{172}}^{-1}$ in  (\ref{m1.11}) and by
using the fact that $\mbox{\ding{172}}^{0}=1$ (see (\ref{3.2.1}))
we obtain
 \beq
f(y+\mbox{\ding{172}}^{-1})=f(y)\mbox{\ding{172}}^{0}+f'(y)\mbox{\ding{172}}^{-1}+
\frac{f''(y)}{2}\mbox{\ding{172}}^{-2}+ \ldots +
\frac{f^{(k)}(y)}{k!}\mbox{\ding{172}}^{-k} + \ldots \label{m1.2}
       \eeq
The uniqueness of the Taylor expansion   allows us to obtain
(\ref{m1}) by equating the first $k+1$ coefficients of \ding{172}
in (\ref{m1.2})  with   grossdigits $c_{0}, c_{-1},  c_{-2},
\ldots c_{-(k-1)}, c_{-k}$ in  (\ref{m1.1}) completing so the
proof. \hfill $\Box$

Let us comment upon the theorem. It describes a situation where a
user needs to evaluate $f(x)$ and its derivatives at a point $x=y$
but analytic expressions of $f(x),$ $f'(x), f''(x), \ldots $ $
f^{(k)}(x)$ are unknown and computer procedures for calculating
$f'(x), f''(x), \ldots f^{(k)}(x)$ are unavailable. Moreover, the
internal structure of the procedure $P$ calculating $f(x)$ can
also be unknown to the user. In this situation, instead of the
usage of, for instance,  traditional formulae (\ref{m0}),
(\ref{m0.1}) for an approximation of $f'(x)$, the user evaluates
$f(x)$ at the point $x=y+\mbox{\ding{172}}^{-1}$ at the Infinity
Computer. Note that if $P$ has been written by the programmer for
the Infinity Computer, then the user just runs $P$ without any
intervention on the code of $P$. In the case when $P$ has been
written for   traditional computers, in order to transfer it to
the Infinity Computer,   variables and constants used in $P$
should be just redeclared as grossnumbers (\ref{3.13}).
Traditional arithmetic operations are then overloaded due to the
rules defined in \cite{Sergeyev_patent,informatica}.

The operation of evaluation of $f(x)$ at the point
$x=y+\mbox{\ding{172}}^{-1}$ returns a number in the form
(\ref{m1}) from where the user can easily obtain   values  of
$f(y)$ and $f'(y), f''(y), \ldots f^{(k)}(y)$ as shown in
(\ref{m2}) without any knowledge of the Taylor expansion of $f(x)$
and of the analytic formulae and computer procedures for
evaluating derivatives. Due to the fact that the Infinity Computer
is able to work with infinite and infinitesimal numbers
numerically, the values $f'(y), \ldots f^{(k)}(y)$ are calculated
exactly at the point $x=y$ without introduction of dangerous
operations (\ref{m0}), (\ref{m0.1}) (or (\ref{m0.2})) related to
the necessity to use finite values of $h$ when one works with
traditional computers. We emphasize also that the user obtains the
function value and the values of the derivatives after calculation
of $f(x)$ at a single point.

It is worthy to notice that numerical operations that the Infinity
Computer performs   when it executes the program $f(x)$ can be
viewed as an automatic rewriting  of $f(x)$ from the basis in $x$
into the basis in $\mbox{\ding{172}}$ by setting
$x=y+\mbox{\ding{172}}^{-1}$ with $y$ being a finite number. The
numerical finite value of $y$ is then combined with other finite
numbers present in the program and they all are collected as
finite coefficients (i.e., grossdigits) of grosspowers of
$\mbox{\ding{172}}$. In some sense this is   similar to
rearrangements that often are executed when one works with
wavelets (see \cite{Walnut}) or  with formal power series (see
\cite{Wilf}).

Let us consider some numerical examples. Their   results   can be
checked by the reader directly on systems using symbolic
calculations (e.g., MAPLE) by taking instead of
$\mbox{\ding{172}}^{-1}$ a symbolic parameter, let say, $a$,
thinking about $a$ as an infinitesimal number and by calculating
then $f(y+a)$ where $y$ is a number. The crucial difference of the
Infinity Computer with respect to  systems executing symbolic
computations consists of the fact that the Infinity Computer works
with infinite, finite, and infinitesimal numbers
\textit{numerically}, not symbolically. Naturally, this feature of
the Infinity Computer becomes very advantageous when one should
execute complex numerical
 computations.

\begin{example}
\label{e_m1}   Suppose that we have a computer procedure
implementing the following function $g(x)=x^3$ as $f(x)=x\cdot
x\cdot x$ and we want to evaluate the values $f(y), f'(y),
f''(y),$ and $f^{(3)}(y)$ at the point $y=5$. The Infinity
Computer executes   the following operations
 \[
f(5+\mbox{\ding{172}}^{-1})=5\mbox{\ding{172}}^{0}1\mbox{\ding{172}}^{-1}\cdot5\mbox{\ding{172}}^{0}1\mbox{\ding{172}}^{-1}\cdot5\mbox{\ding{172}}^{0}1\mbox{\ding{172}}^{-1}=
\]
 \beq
25\mbox{\ding{172}}^{0}10\mbox{\ding{172}}^{-1}1\mbox{\ding{172}}^{-2}\cdot5\mbox{\ding{172}}^{0}1\mbox{\ding{172}}^{-1}=
125\mbox{\ding{172}}^{0}75\mbox{\ding{172}}^{-1}15\mbox{\ding{172}}^{-2}1\mbox{\ding{172}}^{-3}.
\label{m2.3}
       \eeq
From (\ref{m2.3}), by applying (\ref{m2}) we obtain that
\[
f(5)=125,\hspace{3mm}f'(5)=75,\hspace{3mm} f''(5)=2!\cdot 15=30,
\hspace{3mm} f^{(3)}(5)= 3! \cdot 1=6,
\]
that are correct values of $f(x)$ and the derivatives at the point
$y=5$.

Let us check this numerical result analytically by taking a
generic  point $y$. Then we obtain
 \beq
f(y+\mbox{\ding{172}}^{-1})=(y+\mbox{\ding{172}}^{-1})^3
=(y+\mbox{\ding{172}}^{-1})\cdot (y+\mbox{\ding{172}}^{-1}) \cdot
(y+\mbox{\ding{172}}^{-1}) = \label{m2.1}
       \eeq
 \beq
y^3+3y^2\mbox{\ding{172}}^{-1}+3y\mbox{\ding{172}}^{-2}+\mbox{\ding{172}}^{-3}=
y^3\mbox{\ding{172}}^{0}3y^2\mbox{\ding{172}}^{-1}3y\mbox{\ding{172}}^{-2}1\mbox{\ding{172}}^{-3}
. \label{m2.2}
       \eeq
By applying (\ref{m2}) we have the required values
\[
f(y)=y^3,\hspace{3mm}f'(y)=3y^2,\hspace{3mm} f''(y)=2!\cdot 3y=6y,
\hspace{3mm} f^{(3)}(y)= 3! \cdot 1=6.
\]
That coincide with the respective analytical derivatives
calculated at the point $x=y$:
\[
\hspace{25mm} f'(x)=3x^2,\hspace{3mm} f''(x)=6x, \hspace{3mm}
f^{(3)}(x)= 6. \hspace{35mm} \Box
\]

   \end{example}

\begin{example}
\label{e_m2} Suppose  that we have the following function
$g(x)=x+\sin(x)$ and it is represented in the Infinity Computer as
 \beq
  f(x) = x+\widetilde{\sin}(x),
   \label{m3}
   \eeq
where $\widetilde{\sin}(x)$ is a computer implementation of
$\sin(x)$. If we want to evaluate $f(x), f'(x), f''(x),$ and
$f^{(3)}(x)$   at a point $y$, by taking $k=3$ in (\ref{m1}) we
obtain
\[
f(y+\mbox{\ding{172}}^{-1})=
(y+\widetilde{\sin}(y))\mbox{\ding{172}}^{0}(1+\widetilde{\sin}'(y))\mbox{\ding{172}}^{-1}\frac{\widetilde{\sin}''(y)}{2}\mbox{\ding{172}}^{-2}\frac{\widetilde{\sin}^{(3)}(y)}{3!}\mbox{\ding{172}}^{-3},
\]
where the result depends on the way of implementation of
$\widetilde{\sin}(x)$. For example, suppose for the illustrative
purpose that in the neighborhood of the point $y=0$ the Infinity
Computer uses the following simple implementation
 \[
  \widetilde{\sin} (x) = x-\frac{x\cdot x \cdot x}{6}
   \]
being the first two items in the corresponding Taylor expansion.
Then the computer program $f(x)$ becomes
\[
f(x) = x+x-\frac{x\cdot x \cdot x}{6}
\]
and the Infinity Computer with $y=0$ works as follows
\[
f(0+\mbox{\ding{172}}^{-1})=0+\mbox{\ding{172}}^{-1}+0+\mbox{\ding{172}}^{-1}
- \frac{(0+\mbox{\ding{172}}^{-1})\cdot (0+\mbox{\ding{172}}^{-1})
\cdot (0+\mbox{\ding{172}}^{-1}) }{6} =
\]
\[
2\mbox{\ding{172}}^{-1} -
\frac{\mbox{\ding{172}}^{-3}}{6}=2\mbox{\ding{172}}^{-1}
(\mbox{\small-}0.166667)\mbox{\ding{172}}^{-3}.
\]
By applying (\ref{m2}) we have the required values
\[
f(0)=0,\hspace{3mm}f'(0)=2,\hspace{3mm} f''(0)=2!\cdot 0=0,
\hspace{3mm} f^{(3)}(0)= 3! \cdot (\mbox{\small-}0.166667)=-1.
\]
That, obviously, coincide with the respective analytical
derivatives (that, we emphasize this fact again, were not used by
the Infinity Computer)
\[
f'(x)=2-0.5x^2,\hspace{3mm} f''(x)=-x, \hspace{3mm} f^{(3)}(x)= -1
\]
calculated at the point $y=0$.  \hfill $\Box$
  \end{example}
  \begin{example} \label{e_m3} Suppose that we have
a computer procedure $f(x)=\frac{x \cdot x+1}{x}$ implementing the
function $g(x)=\frac{x^2+1}{x}$ and we want to calculate the
values $f(y), f'(y),$ $  f''(y),$ and $f^{(3)}(y)$ at a point
$y=3$. We consider the Infinity Computer that returns grossdigits
 corresponding   to the exponents of grossone from
0 to -3. Then   we have
\[
f(3+\mbox{\ding{172}}^{-1})=\frac{(3+\mbox{\ding{172}}^{-1})\cdot
(3+\mbox{\ding{172}}^{-1}) +1}{3+\mbox{\ding{172}}^{-1}}=
\frac{10\mbox{\ding{172}}^{0}6\mbox{\ding{172}}^{-1}1\mbox{\ding{172}}^{-2}}{3\mbox{\ding{172}}^{0}1\mbox{\ding{172}}^{-1}}=
\]
\[
3.333333\mbox{\ding{172}}^{0}0.888889\mbox{\ding{172}}^{-1}0.037037\mbox{\ding{172}}^{-2}-0.0123457\mbox{\ding{172}}^{-3}.
\]
By applying (\ref{m2}) we obtain that
\[
f(3)=3.333333,\hspace{3mm}f'(3)=0.888889,
\]
\[
f''(3)=2!\cdot 0.037037=0.074074, \hspace{3mm} f^{(3)}(3)= 3!
\cdot (-0.0123457)=-0.074074,
\]
that are   values which one obtains by using explicit analytic
formulae
\[
f(x)=\frac{x^2+1}{x},\hspace{3mm} f'(x)=1-x^{-2},\hspace{3mm}
f''(x)=2x^{-3}, \hspace{3mm} f^{(3)}(x)= -6x^{-4}
\]
for $f(x)$ and its derivatives at the point $x=3$.  \hfill  $\Box$
   \end{example}

%\markboth{Bibliography}{Bibliography}
\bibliographystyle{plain}
\bibliography{XBib_Magic}
%\end{article}
\end{document}